\documentclass{article}
\usepackage[namelimits,intlimits,sumlimits,leqno]{amsmath}
\usepackage{amsfonts}
\usepackage{amssymb}
\usepackage{graphicx}
\newtheorem{theorem}{Theorem}

\newtheorem{proposition}{Proposition}
\newtheorem{lemma}{Lemma}
\newtheorem{definition}{Definition}

\newtheorem{conjecture}{Conjecture}
\newenvironment{proof}{{\par\noindent\bf Proof. }}{$\Box$}

\newcounter{mili}
\newenvironment{milista}{\begin{list}{\roman{mili})}
{\usecounter{mili}}}{\end{list}}
\newcounter{milib}

\makeatletter
\@addtoreset {equation}{section}
\makeatother
\def\leq{\leqslant}
\def\geq{\geqslant}

\def\R{\mathbb{R}}
\def\N{\mathbb{N}}
\def\R{{\mathop{{\rm I}\kern-.2em{\rm R}}\nolimits}}
\def\supp{{\mbox{\rm supp}}}

\newcommand{\negrita}[1]{\mbox{\boldmath $ #1 $}}

\begin{document}

\title{Latin Squares, Partial Latin Squares and 
its Generalized Quotients}
\author{L. Yu. Glebsky and Carlos J. Rubio}
\date{26.03.2003}
\maketitle

\begin{abstract}
A (partial) Latin square is a table of multiplication of a 
(partial) quasigroup.
Multiplication of a (partial) quasigroup may be considered as 
a set of triples.
We give a necessary and sufficient condition when a set of
 triples is a 
quotient of a (partial) Latin square.  
\end{abstract}





\section{Introduction}\label{sec:1}
 
Generalized quotient of quasigroup (quotient with respect to
an equivalence relation which is not a congruence) was 
introduced and studied by A.J.W. Hilton. Such a quotient
is neither a quasigroup nor an algebraic system. 
It may be thought as a multivalued algebraical system or
a set of triples. Theorem~\ref{th.Hilton} (A.J.W. Hilton) 
gives necessary and
sufficient conditions when a set of triples is
a generalized quotient of a quasigroup. Here we
extend it to generalized quotients of partial 
quasigroups. Investigation of generalized quotients leads
to a more general objects -- 3-indexed matrices.
The authors of \cite{G1} found useful the notion of generalized quotient
in their investigations of approximations of algebraic systems.

The article is organized as follows.
In Section~\ref{sec:2} we formulate main results
(Theorem~\ref{th.main1} and Theorem~\ref{th.main2})
about 3-indexed matrices. In Section~\ref{sec:2.5}
we define generalized quotient partial quasigroup (GQPQ)
and generalized uniformly quotient partial quasigroup
(GUQPQ), interpret results about matrices on this language,
and discuss some connection with the theory of hypergraphs.
We give example of a GQPQ which is not a GUQPQ and formulate 
conjectures. In section~\ref{sec:3} we prove results on 
2-indexed matrices, which are used in the proof of 
Theorem~\ref{th.main1}. Section~\ref{sec:4} is devoted
to the proof of Theorem~\ref{th.main1}.




\section{Formulation of the main results}\label{sec:2}

We will deal with  2- and 3-indexed matrices.  
For positive 
integers 
$\displaystyle{n_1,\ldots,n_k}$,
a $k$-indexed $n_1\times\cdots\times n_k$-matrix $M$ 
is a function 
$M:(n_1)\times \cdots \times (n_k)\to \R$,
where $(n)=\{1,\ldots ,n\}$. 
Through this text we will use notation 
$M(X)=\sum_{x\in X} M(x)$, where  
$X\subseteq (n_1)\times \cdots \times (n_k)$ .

We will denote by $T(n_1,\ldots ,n_k)$
the set of all $k$-indexed 
$n_1\times\cdots\times n_k$-matrices
with entries being nonnegative integers.  

We call an $n_1\times\cdots\times n_k$-line   
any set 
$l\subset (n_1)\times \cdots \times (n_k)$ 
such that in all  $k$-tuples in $l$, $k-1$ indexes are fixed 
and the other,
say the $i$-th index, runs over all $(n_i)$.
A line of a $n_1\times\cdots\times n_k$-matrix is the restriction
of this matrix on a $n_1\times\cdots\times n_k$-line. 
If $l$ is a line of $M$, 
$M(l)$ is its line sum.

For $n_1\times n_2$-lines we will use the following names and 
notations.
$$
l^1_{a}=\{(x,a):x\in (n_1)\}\;\; \mbox{(column)}
$$
$$
l^2_{a}=\{(a,x):x\in (n_2)\}\;\; \mbox{(row)}
$$
Similarly, for $n_1\times n_2\times n_3$-lines  
we have
$$
l^1_{ab}=\{(x,a,b):x\in (n_1)\}\;\; \mbox{(horizontal line)}
$$
$$
l^2_{ab}=\{(a,x,b):x\in (n_2)\}\;\; \mbox{(transversal line)}
$$
$$
l^3_{ab}=\{(a,b,x):x\in (n_3)\}\;\; \mbox{(vertical line)}
$$

For a function $f:(n)\to(n)$, {\em the graph of $f$} is 
the matrix 
$\Gamma\in T(n,n)$ such that 
$$
\Gamma(i,j)=\left\{\begin{array}{lll}
                           1, & \mbox{if} & f(i)=j,\\
                           0, & \mbox{if} & f(i)\neq j. 
               \end{array}\right. 
$$

It is easy to see that the following proposition holds.

\begin{proposition}\label{prop.simpl2}
A matrix $M\in T(n,n)$ is the graph of a permutation if and 
only if every line sum of $M$ equals $1$. We will call
such a matrix to be a permutation one.
\end{proposition}

An analogue of this proposition for 3-indexed matrices leads 
to quasigroups and Latin squares.

\begin{definition}\label{def.2.1}
A quasigroup $(Q,\star)$ is an algebraic system
$Q$ with a binary operation $\star$ such that
\begin{milista}
\item equation $x\star a=b$ has a unique solution with 
respect to 
$x$ for
all $a,b\in Q$,
\item equation $a\star x=b$ has a unique solution with respect 
to
 $x$ for
all $a,b\in Q$.
\end{milista}
\end{definition}

This definition implies immediately the following

\begin{proposition}\label{prop.simpl3}
A matrix $M\in T(n,n,n)$ is the graph of a quasigroup operation
$\star$ on $(n)$ ($M(i,j,k)=1$, if $i\star j=k$ and $M(i,j,k)=0$ 
otherwise) 
if and 
only if every line sum of $M$ equals $1$. We will call
such a matrix to be a Latin square.
\end{proposition}

For 2-indexed matrices the following lemma is well-known e.g.\cite{Ryser}.
\begin{lemma} \label{lemma.konig}
Let $M\in T(n,n)$. Let each line sum of $M$ equals $k$, 
where $k>0$. Then
$$
\supp(M)\supseteq\supp(P)
$$
for a permutation matrix $P$.
\end{lemma}
This Lemma easily implies
\begin{lemma} \label{lemma.3.3}
Let $M\in T(n,n)$. Let each line sum of $M$ equals $k$, where $k>0$. 
Then
$$M=P_{1} + P_{2} + \cdots + P_{k},$$
where each $P_{i}$ is a permutation matrix.
\end{lemma}

One may formulate the following ``generalization'' of 
Lemma~\ref{lemma.konig}\\

(A) {\em Let $M\in T(n,n,n)$. Let each line sum of $M$ equals $k$,
 whith $k>0$. Then
$$
\supp(M)\supseteq\supp(L),
$$
for a Latin square $L$.}\\

Statement (A) is not true. Indeed, consider the
following $3\times 3\times 3$ matrix
\begin{center}
\includegraphics[width=2.5in]{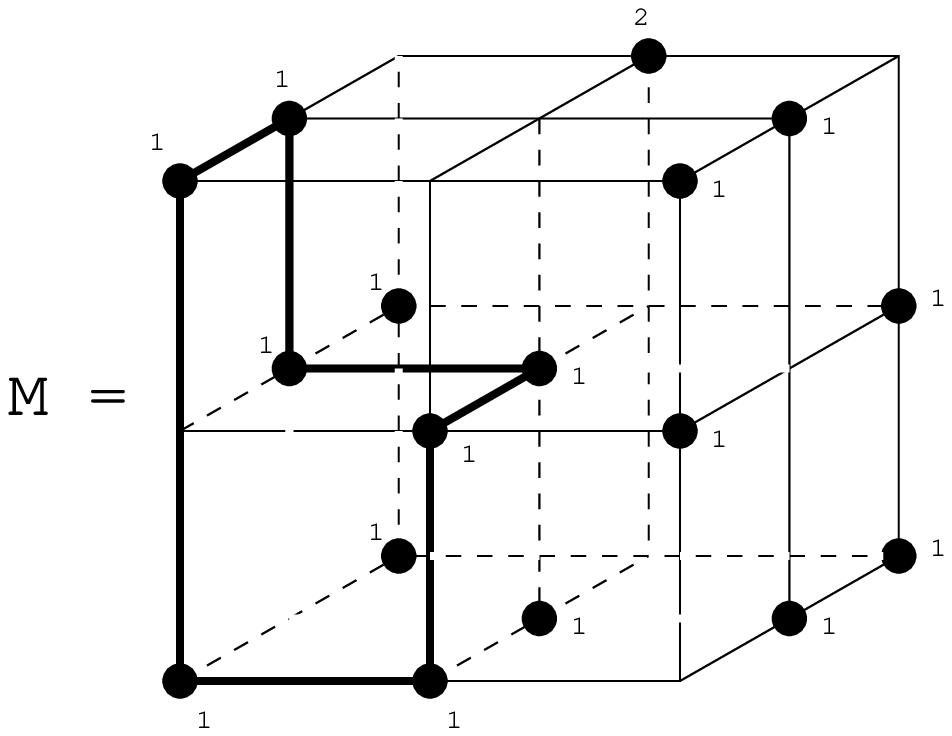}
\end{center}
where every line sum of $M$ equals 2.
The existence of the odd cycle ${\cal C}$ in $M$ 
(marked bold) 
implies that $\supp(M)\not\supseteq\supp(L)$  
for any Latin square $L$. Indeed, 
let $\supp(M)\supseteq \supp(L)$ for some Latin square $L$. 
Then $\supp(L)$ has
to contain only  one dot marked of every line of $\cal C$. But 
this is impossible because ${\cal C}$ is an odd cycle.\\ 
{\bf Remark.} {\em A set of triples, as a hypergraph may not have 
a (dual) K\"onig
property. On the contrary, any set of pairs is a balanced 
hypergraph
and satisfies a (dual) K\"onig property, i.e. $\rho=\bar\alpha$, see
Section~\ref{sec:2.5} and \cite{Berge} }\\
Nevertheless, there is some connection of matrices described in
statement (A) with Latin squares through quotients.

Let $M:(n_1)\times \cdots \times (n_k)
\rightarrow \R$ be a 
$k$-indexed matrix and $\sigma = \{P_1, \ldots ,P_r\}$ 
be an (ordered) partition of $(n_i)$. We define 
the quotient matrix
$$
M{\circ}_i \sigma: 
(n_1)\times \cdots \times (n_{i-1})\times (r)\times 
(n_{i+1})\times \cdots \times (n_k)\rightarrow \R,
$$
by the formula
$$
M{\circ}_i \sigma (x_1,\ldots, x_i, \ldots, x_k) = 
M(\{x_1\}\times \cdots \times P_{x_i}\times \cdots \{x_k\}).
$$\\

Example. Let
$$
M= \left( \begin{array}{cccc}
0 & 3 & 3 & 1\\
5 & 2 & 4 & 0\\
1 & 1 & 0 & 1\\
2 & 3 & 5 & 0
\end{array} \right)
$$
and $\sigma=\{\{1,2\},\{3,4\}\}$. Then
$$
M\circ_1\sigma = \left( \begin{array}{cccc}
5 & 5 & 7& 1\\
3 & 4 & 5 & 1
\end{array} \right),\;\;\;
    M\circ_2\sigma =
\left( \begin{array}{cc}
3 & 4\\
7 & 4\\
2 & 1\\
5 & 5 
\end{array} \right)
$$
and
$$
  (M\circ_1\sigma)\circ_2\sigma=(M\circ_2\sigma)\circ_1\sigma=
\left( \begin{array}{cc}
10 & 8\\
7 & 6 
\end{array} \right).
$$

Let $L\in T(n,n,n)$ be a Latin square and 
$\sigma=\{P_1,P_2,\ldots,P_k\}$ be a partition of $(n)$,
let $M=((L\circ_1\sigma)\circ_2\sigma)\circ_3\sigma\in T(k,k,k)$.
Then it is easy to check that 
$$
M(l^t_{ij})=|P_i|\cdot |P_j|,
$$
for every $k\times k\times k$-line $l^t_{ij}$, $t=1,2,3$,
$i,j\in(k)$.   
It was proved by  Hilton \cite{Hilton1,Hilton2} that the inverse
statement is also true.

\begin{theorem}[A.J.W.Hilton] 
\label{th.Hilton}
Let $M\in T(k,k,k)$ and $r_1, r_2,\ldots, r_k$ be 
positive integers such that 
$$
M(l^t_{ij})=r_ir_j
$$
for $t=1,2,3$ and $i,j\in (k)$.
Then $M=((L\circ_1\sigma)\circ_2\sigma)\circ_3\sigma$ for
a Latin square $L\in T(n,n,n)$ and a partition 
$\sigma=\{P_1,\ldots,P_k\}$  of $(n)$ such that 
 $|P_i|=r_i$ and $n=\sum_ir_i$.
\end{theorem}

In the paper we generalize this theorem for partial 
Latin squares (and partial quasigroups).

\begin{definition}\label{defi.2.2}
Let $Q$ be a finite set and $S\subseteq Q\times Q$. 
A partial $S$-quasigroup on $Q$ is a partial binary operation 
$\star$ on $Q$ such that
\begin{milista}
\item $S\subseteq Dom(\star)$.
\item equation $x\star a= b$ has at most one solution with 
respect to $x$ for all $a,b\in Q$.
\item equation  $a\star x =b$ has at most one solution 
with respect to $x$ for all $a,b\in Q$.
\end{milista}
\end{definition}

\begin{proposition}\label{prop.simpl3.s}
A matrix $M\in T(n,n,n)$ is the graph of a partial $S$-quasigroup 
operation on (n) if and 
only if every line sum of $M$ is no more than 1 and
$$
M(l^3_{ij})=1
$$
for every $(i,j)\in S$.
We will call
such a matrix to be a partial $S$-Latin square.
\end{proposition}

Let $L\in T(n,n,n)$ be a partial $S$-Latin square, 
$\sigma=\{P_1,\ldots,P_k\}$ be a partition of $(n)$,
and $M=((L\circ_1\sigma)\circ_2\sigma)\circ_3\sigma\in T(k,k,k)$.
Then it is easy to verify that  
$$
M(l^t_{ij})\leq |P_i|\cdot |P_j|
$$
for every $k\times k\times k$-line 
$l^t_{ij}$, $t=1,2,3$,
$i,j\in(k)$  and
$$
M(l^3_{ij})= |P_i|\cdot |P_j|,
$$
if $P_i\times P_j\subseteq S$.

\begin{theorem} \label{th.main1}
Let $M\in T(k,k,k)$, $S\subseteq(k)\times(k)$
 and $r_1, \ldots, r_k$ be 
positive integers such that 
$$
M(l^t_{ij})\leq r_ir_j
$$
for $t=1,2,3$, and
$$
M(l^3_{ij})= r_ir_j,
$$
for $(i,j)\in S$.\\
Then $M=((L\circ_1\sigma)\circ_2\sigma)\circ_3\sigma$ for
a partial $S'$-Latin square $L\in T(n,n,n)$, a partition 
$\sigma=\{P_1,\ldots,P_k\}$  of $(n)$ such that 
 $|P_i|=r_i$, $n=\sum_ir_i$, and 
$$
S'=\bigcup_{(i,j)\in S} P_i\times P_j.
$$
\end{theorem}

\noindent If we substitute $S=(k)\times(k)$ in 
Theorem~\ref{th.main1} 
we get Theorem~\ref{th.Hilton}. (If for all vertical lines
one has equalities then one has equalities for all lines.)
The uniform partial case of Theorem~\ref{th.main1}, where
$r_1=\cdots=r_k=r$, may be generalized for real-valued
matrices.

\begin{theorem} \label{th.main2}
Let $M:(k)\times(k)\times (k) \to \R^+$, 
$\beta\in \R^+$ and
$S\subseteq(k)\times(k)$
 such that
$$
M(l^t_{ij})\leq \beta
$$
for $t=1,2,3$ and $i,j\in (k)$, and
$$
M(l^3_{ij})= \beta,
$$
for $(i,j)\in S$.\\
Then $\supp(M)=
\supp(((L\circ_1\sigma)\circ_2\sigma)\circ_3\sigma)$ for
a partial $S'$-Latin square $L\in T(n,n,n)$, a partition 
$\sigma=\{P_1,\ldots,P_k\}$ of $(n)$ such that $|P_i|=|P_j|$ for every 
$i\neq j$, and
$$
S'=\bigcup_{(i,j)\in S} P_i\times P_j.
$$
\end{theorem}
\begin{proof}
Let $M$ and $\beta$ satisfy the conditions of the theorem. 
We can write equalities and strict inequalities separately.
Consider non-zero elements of $M$ and $\beta$ as variables.
Then this system of equalities and (strict) inequalities has 
a rational solution. Multiplying this solution by a proper
integer, we construct a matrix $M'\in T(k,k,k)$, 
$\supp(M')=\supp(M)$, satisfying the 
conditions of Theorem~\ref{th.main1}.
\end{proof}

As we see, the crucial step in
the proof of Theorem \ref{th.main2} is to show the existence of a 
rational solution.
For the non-uniform case these equations will be 
nonlinear (quadratic). So, general consideration cannot
prove that the existence of a real solution implies the 
existence of a rational one. We don't know so far if a 
non-uniform version of 
Theorem~\ref{th.main2} is valid.



\section{Generalized quotient quasigroup} \label{sec:2.5}

Let $Q$ be a finite set and $\sigma$ an equivalence relation
on $Q$ which we will identify with the partition of $Q$ 
by equivalence classes. So,
 $\sigma=\{Q_1,Q_2,...,Q_k\}$. Let $X\subseteq Q^r$. Define
weak ($X/^w\sigma\subseteq \{1,2,...k\}^r$) 
and strong ($X/^s\sigma\subseteq\{1,2,...k\}^r$) quotient of $X$:
$$
X/^w\sigma=\{<i_1,i_2,...,i_k>\; :\; 
Q_{i_1}\times Q_{i_2}\times\cdots\times Q_{i_k}\cap X
\neq\emptyset\}
$$
$$
X/^s\sigma=\{<i_1,i_2,...,i_k>\; :\; 
Q_{i_1}\times Q_{i_2}\times\cdots\times Q_{i_k}\subseteq X \}
$$
For example, if $\star\subseteq Q^3$ is a quasigroup operation 
on $Q$ and $\sigma$ -- a congruence relation 
(i.e. it preserves the operation $\star$) then 
$\star/^w\sigma=\star/^s\sigma=\star/\sigma$ is a quotient 
quasigroup operation.
\begin{definition}
Let $\sigma=\{Q_1,Q_2, ...,Q_k\}$ be an equivalence relation on 
$Q$. We will call $\sigma$ to be uniform iff all $Q_i$ have
the same cardinality.

Let $(Q,\star)$ be a quasigroup. A set 
$\star/^w\sigma$ will be called
a generalized quotient  quasigroup (GQQ). 
For uniform $\sigma$ a set 
$\star/^w\sigma$ will be called
a generalized uniformly quotient  quasigroup (GUQQ).   

Let $(Q,\star)$ be an $S$-quasigroup ($S\subseteq Q^2$) and
$\sigma$ be an equivalence relation on $Q$ 
(not necessarily a congruence).
A set $\star/^w\sigma$ will be called
a $S/^s\sigma$-generalized quotient partial quasigroup 
($S/^s\sigma$-GQPQ) or a $S/^s\sigma$-generalized 
uniform quotient partial quasigroup ($S/^s\sigma$-GQUPQ) in the 
case of uniform $\sigma$.  
\end{definition}

Theorem~\ref{th.main1} and Theorem~\ref{th.main2} 
have the following obvious 
interpretation:\\
{\it 
 $H\subseteq (k)^3$ 
is a $S$-GQPQ ($S$-GUQPQ) if and only if there exists a matrix $M$, 
$\supp(M)=H$, satisfying conditions of Theorem~\ref{th.main1}
(Theorem~\ref{th.main2}).
} 

Now we give an interpretation of our results on the language of
hypergraphs. 
A set of triples $H\subseteq (k)^3$ has a natural structure of 
a hypergraph if we consider lines as  edges. Precisely,
with $H$ we associate hypergraph $(V,E)$ with $V=H$ and 
$E=\{l\cap H\;:\;l-\mbox{line}\}$. Several useful
numeric characteristic of hypergraphs are known. 
We are interested in 3
of them: the covering number $\rho$, the independent number
$\bar\alpha$ and  the fractional independent number $\alpha^*$,
for the general definition, see \cite{Berge}.
For the case of a set of triples $H\subseteq (k)^3$ 
these numbers have the 
following meaning:\\
$\rho(H)$ is the minimum number of lines, covering $H$, i.e. 
$\rho(H)=\min\{|R|\;:\; R\;\mbox{is a set of lines, and}\; 
H\subseteq\cup R\}$;\\
$\bar\alpha(H)=\max\{|X|\;:\;X\subseteq H\;
\mbox{and}\;|X\cap l|\leq 1\;\mbox{for every line}\;l\}$\\
$\alpha^*(H)=\max\{M(H)\;:\;M:(k)\times (k)\times (k)\to\R^+,\; 
\supp(M)\subseteq H\;\mbox{and}\;M(l)\leq 1 
\;\mbox{for any line}\; l\}$.\\
From the theory of hypergraphs, see \cite{Berge}, it follows
that $\bar\alpha(H)\leq\alpha^*(H)\leq\rho(H)$. One immediately
verifies the following 
\begin{proposition}
$H\subseteq (k)^3$ contains an 
$S$-quasigroup if and only if
$\bar\alpha(H\cap (S\times (k)))=|S|$.
\end{proposition}
Without loss of generality one
can put $\beta=1$ in Theorem~\ref{th.main2}. It implies     
\begin{proposition}
$H\subseteq (k)^3$ contains an 
$S$-GUQPQ if and only if\\
$\alpha^*(H\cap (S\times (k)))=|S|$.
\end{proposition}
The following example shows that there exists a GQQ which 
does not contain a GUQQ.
\begin{center}
\includegraphics[width=2.5in]{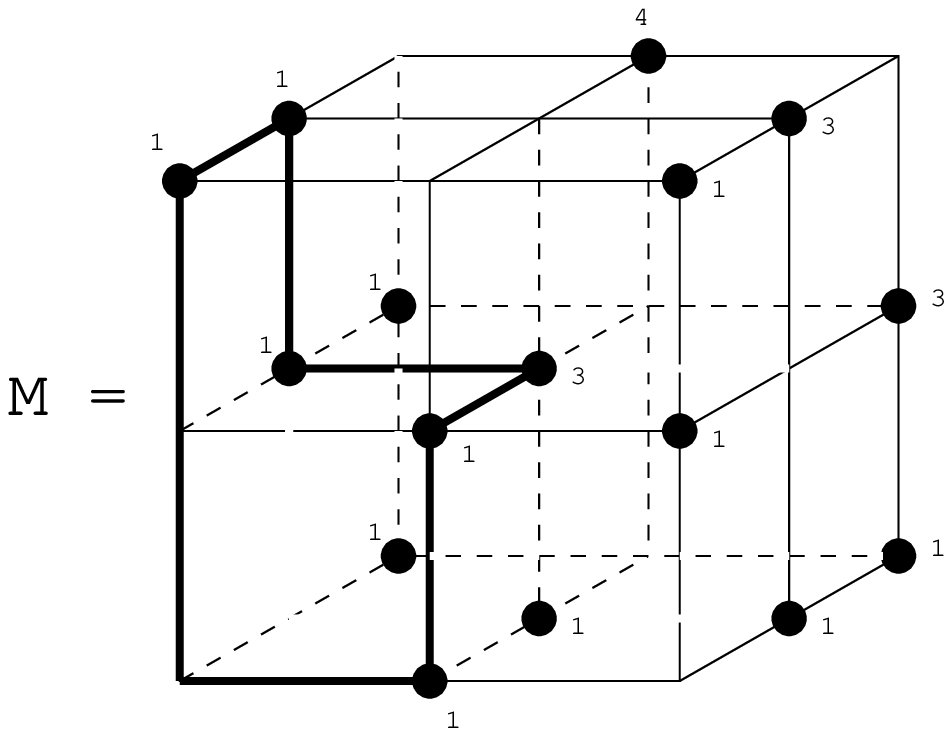}
\end{center}
The numbers show that the set marked by black dots is a GQQ
if we put $r_1=1$, $r_2=2$ and $r_3=2$ in Theorem~\ref{th.Hilton}.
On the other hand if the set were a GUQQ one could put numbers, 
such that sums along every line are the same. To see that
it is impossible one can try to put numbers along the odd
cycle $\cal C$ (marked bold). 

We wonder, if the following conjectures are true.
\begin{conjecture}
If $H\subseteq (k)^3$ contains an 
$S$-GQPQ then
$\rho(H\cap S\times (k))=|S|$.
\end{conjecture}
\begin{conjecture}
If $\rho(H\cap S\times (k))=|S|$ then
$H\subseteq (k)^3$ contains an 
$S$-GQPQ.
\end{conjecture}




\section{Decomposition of 2-indexed matrices of nonnegative 
integers}\label{sec:3}

\noindent Let $M\in T(m,n)$. Let the sum of row $i$ of $M$ be 
denoted by $r_{i}$ and let the sum of column $j$ of $M$ be
denoted by $s_{j}$. We call the vector
$$
R=(r_{1}, \ldots , r_{m})
$$
the {\em row sum vector} and the vector
$$
S=(s_{1}, \ldots , s_{n})
$$
the {\em column sum vector} of $M$.

\noindent The vectors $R$ and $S$ determine the class
$${\cal C}(R,S)$$
consisting of all matrices of size $m$ by $n$, whose entries 
are nonnegative integers, with row sum vector $R$ and column 
sum vector $S$. There is a simple  necessary and sufficient condition under
 which the class ${\cal C}(R,S)$ is nonempty.

We will use the notation 
$\negrita{x} = (x_{1}, \ldots , x_{m})$ for vectors 
such that $x_{i}\in \N$ where $\N = \{0,1,2, \ldots \}$. 
If $\negrita{x}, \negrita{y}\in {\N}^m$, we will say that 
$\negrita{x}\leq \negrita{y}$ if and only if $x_{i}\leq y_{i}$ 
for every $i$.  
Finally we define $\negrita{x}+\negrita{y}=(x_{1}+y_{1}, \ldots , 
x_{m}+y_{m}), k\negrita{x}=(kx_1,\ldots,kx_m)$ and 
$|\negrita{x}|=\sum_{i=1}^m |x_{i}|$ for every 
$\negrita{x}\in {\N}^m$.


\begin{proposition}\label{prop.3.6}
Let $\negrita{x}\in {\N}^m$ and let $N\in \N$ such that 
$|\negrita{x}|\geq N$. Then there exists $\negrita{y}\in {\N}^m$ 
such that $\negrita{y}\leq \negrita{x}$ and $|\negrita{y}|=N$.
\end{proposition}
\begin{proof}
We will use induction on $N$. 
For $N=1$ the proposition is trivial. If $N=N_{1}+N_{2}$, 
we can find $\negrita{y}_{1}\in {\N}^m$ such that 
$\negrita{y}_{1}\leq \negrita{x}$ and  $|\negrita{y}_{1}|=N_{1}$, 
and $\negrita{y}_{2}\in {\N}^m$ such that 
$\negrita{y}_{2}\leq \negrita{x} - \negrita{y}_{1}$ and 
$|\negrita{y}_{2}|=N_{2}$. Making 
$\negrita{y}=\negrita{y}_{1}+\negrita{y}_{2}$ the proposition 
follows. 
\end{proof}

\begin{proposition}\label{prop.3.7}
Let $\negrita{x}\in {\N}^m$ and suppose that $|\negrita{x}|=rk$ 
for some positive integers $r,k$. Then there exist 
$\negrita{y}_{1}, \negrita{y}_{2}, \ldots , 
\negrita{y}_{r}\in {\N}^m$ such that $|\negrita{y}_{i}|=k$ and 
$\negrita{x}=\negrita{y}_{1}+\negrita{y}_{2}+\cdots +
\negrita{y}_{r}$.
\end{proposition}
\begin{proof}
The proof follows from
Proposition~\ref{prop.3.6}.  
\end{proof}

\begin{proposition} \label{prop.3.8}
Let $R=(r_1,r_2, \ldots , r_m)$ and 
$S=(s_1,s_2, \ldots , s_n)$ be nonnegative vectors with integer 
entries. Then the class ${\cal C}(R,S)$ is nonempty if and only 
if 
$$
|R|=|S|.
$$
\end{proposition}
\begin{proof}
$(\Rightarrow)$: Let the class ${\cal C}(R,S)$ contain a matrix 
$M$. Then $|R|$ and $|S|$ are 
both equal to the sum of all the entries of $M$, 
and therefore $|R| = |S|$.

$(\Leftarrow)$: We will use induction on $m$. 
For $m=1$ the proposition is trivial. 
Let $m>1$ and suppose that the proposition is true for $m-1$.
Suppose that $|S|=|R|$. It is clear that $|S|\geq r_1$. 
Then by Proposition~\ref{prop.3.6}, there exists 
$\negrita{x}\in \N^n$ such that $\negrita{x}\leq S$ and 
$|\negrita{x}|=r_1$. Then
$$
M=\left(\begin{array}{c}
\negrita{x}\\
M'
\end{array}\right)
$$
where $M'\in {\cal C}(R',S'), R'=(r_2,\ldots,r_m)$ and $S'=S-\negrita{x}$.

Hence $|R'|=|R|-r_1, |S'|=|S|-|\negrita{x}|=|R'|$, and the  matrix 
$M'$ exists by the induction hypothesis.
\end{proof}

\begin{lemma}\label{lemma.3.5}
If $M\in {\cal C}(kR, kS)$, 
 then
$$
M=Q_{1}+Q_{2}+\cdots +Q_{k},
$$
where $Q_{i}\in {\cal C}(R, S)$ for every $i=1,2,\ldots,k$. 
\end{lemma}
\begin{proof}
The idea of the proof is the following: 
Using Proposition~\ref{prop.3.7}, substitute every  row $i$ of 
$M$ by its decomposition in $r_i$ rows in such a way that
the sum of every row equals $k$. Then  substitute
every column $j$ by its decomposition in $s_j$ columns
in such a way that the sum of every column equals $k$.
The resulting matrix satisfies Lemma~\ref{lemma.3.3}. 
Then, using the matrices $P_{i}$ of Lemma~\ref{lemma.3.3}, 
construct 
$Q_{i}$ of Lemma~\ref{lemma.3.5}.

Let $M:(m)\times (n)\rightarrow \N$ 
and $\negrita{m}_{i}=(m_{i1}, \ldots , m_{in})$ be the 
$i$th. row of $M$, $1\leq i\leq n$. 
\noindent Since $|\negrita{m}_{i}|=r_i k$, Proposition~\ref{prop.3.7} 
implies 
that there 
exist 
$\negrita{m}^1_{i}, \negrita{m}^2_{i}, \ldots , 
\negrita{m}^{r_i}_{i}\in {\N}^n$ such that 
$|\negrita{m}^j_{i}|=k$ for each $1\leq j\leq r_i$, and 
$\negrita{m}^1_{i} + \negrita{m}^2_{i} + \cdots + 
\negrita{m}^{r_i}_{i} = \negrita{m}_{i}$. 
From here, it follows that there exists 
$M':(m')\times (n)\rightarrow \N$ 
such that each row sum of $M'$ equals $k$ and 
$M=M'{\circ}_1 \sigma$, where $\sigma=\{P_1,\ldots,P_m\}$ is a partition of 
$(m')$ and $|P_i|=r_i$. 
Similar considerations show that  
$M'=M''{\circ}_2 \sigma'$, such that each row and column sum of $M''$ 
equals $k$. 
Then by Lemma \ref{lemma.3.3} we have that
\begin{eqnarray}
M''=P_{1}+P_{2}+\cdots +P_{k}.
\end{eqnarray}
It is easy to verify that $Q_i=(P_i\circ_1 \sigma)\circ_2 \sigma'$, 
$i=1,\ldots, k$, 
satisfy the lemma. 
\end{proof}

Let  $R=(r_1,r_2,\cdots,r_m)$, $S=(s_1,s_2,\cdots s_n)$, 
$I\subseteq (n)$. Denote by ${\cal C}_I'(R,S)$ the union of all 
${\cal C}(R',S')$ such that $R'\leq R$, $S'\leq S$ and $s_i=s_i'$ for
$i\in I$. 
\begin{lemma}\label{lemma.4}
Let $M\in {\cal C}_I'(kR,kS)$ such that $|R|=|S|$.
Then
$$M=Q'_1+Q'_2+\cdots+Q'_k,$$
where $Q'_i\in {\cal C}_I'(R,S)$.
\end{lemma}
\begin{proof}
Let $M\in {\cal C}_I'(kR,kS)$ such that $|R|=|S|$. Let $R'$ and $S'$ 
be the row sum and column sum vectors of $M$. Then, $R'\leq kR$,  $S'\leq kS$ 
and from Proposition~\ref{prop.3.8} 
we have that $|R'|=|S'|$. Then $|kR-R'|=|kS-S'|$.
Hence it follows  from Proposition~\ref{prop.3.8} that there exists 
$M'\in {\cal C}(kR-R',kS-S')$. Then $M+M'\in {\cal C}(kR,kS)$, and from 
Lemma~\ref{lemma.3.5} it follows that
$$
M+M' = Q_1+Q_2+\cdots + Q_k,
$$
where $Q_i\in {\cal C}(R,S)$ for every $i=1,\ldots,k$. 
Now, in order to construct $Q_i'$, we have to decrease some entries of 
$Q_i$. 
\end{proof}




\section{Proof of Theorem 2}\label{sec:4}


\noindent Let $M\in T(k,k,k)$, $S\subseteq (k)\times (k)$ and $r_1,\ldots,r_k$ 
be nonnegative integers such that $M(l^t_{ij})\leq r_ir_j$ for $t=1,2,3$, and  
$M(l^3_{ij})=r_ir_j$ for $(i,j)\in S$.
Take $n=r_1+r_2+\cdots r_k$ and partition 
$\sigma=\{P_1,P_2,\cdots, P_k\}$ of $(n)$ such that $|P_i|=r_i$. 
We will consecuently construct 
$M_1\in T(k,k,n), M_2\in T(k,n,n)$ and $M_3\in T(n,n,n)$ such that
\begin{milista}
\item $M=M_1{\circ}_3\sigma,\; M_1=M_2{\circ}_2\sigma,\; 
M_2=M_3{\circ}_1\sigma$;
\item $M_1(l^3_{ij})=r_ir_j$ if $(i,j)\in S$, 
$M_1(l^3_{ij})\leq r_ir_j$, and $M'_1(l^t_{ij})\leq r_i$ for $t=1,2$;
\item $M_2(l^3_{ij})=r_i$ for 
$\displaystyle{(i,j)\in \bigcup_{(i,k)\in S} \{i\}\times P_k}$, 
$M_2(l^3_{ij})\leq r_i$, $M_2(l^2_{ij})\leq r_i$, and $M_2(l^1_{ij})\leq 1$;
\item $M_3(l^t_{ij})\leq 1$ for $t=1,2,3$, and $M_3(l^3_{ij})=1$ 
if $\displaystyle{(i,j)\in\bigcup_{(i,j)\in S} P_i\times P_j}$.
\end{milista}

\noindent {\em Construction of $M_1$.}

\vspace{1ex}

\noindent For every $c$ fixed, $M'_c=M(\cdot,\cdot,c)\in T(k,k)$ such that 
$M'_c(l^t_i)\leq r_ir_c$ for $i\in (k)$ and $t=1,2$. 
So, by Lemma \ref{lemma.4} we can write
$$
M(\cdot,\cdot,c)=Q_{{\alpha}_1}+\cdots+Q_{{\alpha}_{r_c}},
$$
such that $Q_{{\alpha}_m}\in T(k,k)$ and $
Q_{{\alpha}_m}(l^t_i)\leq r_i$ for 
$i\in (k)$ and $t=1,2$. 
One can choose $\alpha_i$ such that 
$\{\alpha_1,\alpha_2,\ldots,\alpha_{r_c}\}=P_c$.
Doing the same for all $c$, we will have matrices 
$Q_1,\ldots,Q_n\in T(k,k)$. Let $M_1(a,b,c)=Q_c(a,b)$.

\noindent {\em Construction of $M_2$.}

\vspace{1ex}

\noindent For every $c$ fixed, $M'_c=M_1(\cdot,c,\cdot)\in T(k,n)$ 
such that $M'_c(l^1_i)\leq r_ir_c$, $M'_c(l^2_i)\leq r_c$, and 
$M_1(l^2_i)=r_c$ 
for $(i,c)\in S$. By Lemma \ref{lemma.4}, we can write
$$
M_1(\cdot,c,\cdot)=Q_{{\alpha}_1} + \cdots + Q_{{\alpha}_{r_c}},
$$
such that $Q_{{\alpha}_m}\in T(k,n)$, 
$Q_{{\alpha}_m}(l^2_i)\leq 1$, $Q_{{\alpha}_m}(l^1_i)\leq r_c$, and 
$Q_{{\alpha}_m}(l^2_i)=1$ if $(i,c)\in S$. 

One can choose $\alpha_i$ such that 
$\{\alpha_1,\alpha_2,\ldots,\alpha_{r_c}\}=P_c$.
Doing the same for all $c$, we will have matrices 
$Q_1,\ldots,Q_n\in T(k,n)$. Let $M_2(a,c,b)=Q_c(a,b)$.

Construction of $M_3$ is similar.
\noindent It is clear that $L=M_3$ satisfies the theorem.






\begin{thebibliography}{}
%
\bibitem{Berge} C.Berege Hypergraphs: Combinatorics of Finite 
Sets. North Holland Mathematical Library, V 45,
Elsevier Science Publishers, 1989.
%
\bibitem{Hilton2} J.K. Dugdale, A.J.W. Hilton, J. Wojciechowski,
Fractional latin squares, simplex algebras, and generalized 
quotients. Journal of statistical planning and inference. 
86, 457-504 (2000).
%
\bibitem{G1} L.Yu. Glebsky, E.I. Gordon, On approximation of topological
groups by finite algebraic systems, preprint  math.GR/0201101,  
http://xxx.lanl.gov/. To be published in Illinois Math. J.
%
%
%
\bibitem{Hilton1} A.J.W. Hilton, Outlines of latin squares. Ann. Discrete 
Math. 34, 225-242 (1987)
%
\bibitem{Ryser} H. J. Ryser, Combinatorial mathematics 
(The Carus Mathematical Monographs, 15) 
The Mathematical Association of America, 1963.
\end{thebibliography}
\end{document}